\documentclass[oneside]{amsart}
\usepackage{amsmath, amssymb,color,amsthm,graphicx,cite}
\usepackage{color,bm}
\newtheorem{theorem}{Theorem}[section]
\newtheorem{proposition}[theorem]{Proposition}
\newtheorem{lemma}[theorem]{Lemma}
\newtheorem{corollary}[theorem]{Corollary}

\theoremstyle{definition}
\newtheorem{definition}{Definition}

\def\Z{\mathbb{Z} }
\def\R{\mathbb{R} }

\def\T{\mathbb{T} }

\def\nbd{neighborhood }

\def\Sv{\mathop{\mathrm{Sing}}(v)}

\def\Pv{\mathop{\mathrm{Per}}(v)}
\def\Cv{\mathop{\mathrm{Cl}}(v)}
\def\F{\mathcal{F} } 
\def\Ct{\mathop{\mathrm{Cl}}(\tau)}

\author{Tomoo Yokoyama}
\date{\today}
\address{Applied Mathematics and Physics Division, Gifu University, Yanagido 1-1, Gifu, 501-1193, Japan\\}
\email{tomoo@gifu-u.ac.jp}
\thanks{The author was partially supported by JSPS Grant Number 20K03583}

\makeatletter
\@namedef{subjclassname@2020}{%
  \textup{2020} Mathematics Subject Classification}
\makeatother


\title[Morse hyper-graphs]{Morse hyper-graphs of topological spaces and decompositions}
\begin{document}
\maketitle

\begin{abstract}
The cell complex structure is one of the most fundamental structures in topology and combinatorics, the Morse decomposition of a dynamical system analyzes the global gradient behavior, and the Reeb graph of a function is an elementary tool in Morse theory to represent the global connection and is also used to analyze continuous and discrete data in topological data analysis. In this paper, we unify three concepts of cell complexes, Morse decompositions, and Reeb graphs into a concept. In fact, we introduce topological invariants for topological spaces and decompositions, which are analogous to abstract (weak) orbit spaces and Morse graphs for flows. To achieve these, we define analogous concepts of recurrence and ``chain-recurrence'' for topological spaces and decompositions such that the original and new recurrences correspond to each other for a flow orbits on a locally compact Hausdorff space and the orbit space. We show that the Morse hyper-graphs exist for any topological spaces and any invariant decompositions (e.g. compact foliated spaces). Moreover, the abstract weak element spaces generalize the Reeb graphs of Morse functions and the Morse (hyper-)graphs and refine abstract cell complexes.
\end{abstract}

\section{Introduction}
The cell complex structure is one of the most fundamental tools in algebraic topology and combinatorics. 
For instance, the CW complex structure, which is an example of the cell complex structure, was introduced by Whitehead \cite{whitehead1949combinatorial} and is widely developed. 
Moreover, the simplicial complex structure, which is also an example of the cell complex structure and is led from the barycentric subdivision of Poincar\'e, is the foundation of the simplicial homology. 
In addition, the existence of triangulations (i.e. the simplicial complex structures) of topological spaces is a classical problem in topology. 
Thus the cell complex structure is a fundamental structure in topology and combinatorics. 
In this paper, we generalize such an abstract cell structure by using concepts of dynamical systems. 

From a dynamical system's point of view, Birkhoff introduced the concepts of recurrent points and described the limit behavior of orbits \cite{Birkhoff}. 
Conley defined a weak form of recurrence, called chain recurrence, for a flow on a compact metric space \cite{Conley1988}.
The Conley theory says that dynamical systems on compact metric spaces can be decomposed into blocks,  each of which is a chain recurrent one or a gradient one.
Then this decomposition is called the Morse decomposition, and implies a directed graph, called a Morse graph, which can capture the gradient behaviors. 
The Morse decompositions for mappings and semiflows are developed for arbitrary metric spaces \cite{franks1988pb,hurley1991chain,hurley1992noncompact,hurley1995chain}. 
Moreover, the Morse decomposition for set-valued dynamical systems is essentially extend by Mcgehee\cite{mcgehee1992attractors} and was developed in several ways \cite{caraballo2015morse,da2016morse,li2007morse,wang2020morse}. 
In addition, the Morse decompositions for semigroups and set-valued semiflows \cite{barros2010attractors,barros2012dynamic,bortolan2013skew,barros2010finest,da2017morse,patrao2007morse,patrao2007semiflows,rybakowski2012homotopy}, random dynamical systems \cite{caraballo2012gradient,crauel2004towards,ju2018strong,lin2021morse,liu2005random,liu2007random,liu2007random3,liu2008attractor}, nonautonomous set-valued dynamical systems \cite{aragao2013non,wang2012morse,wang2014morse}, and combinatorial dynamical systems on simplicial complexes \cite{bogdan2020link,mrozek2017conley,tamal2019pers}, are also introduced and studied in various ways. 
%
In \cite{yokoyama2021refinements}, the Morse graph of dynamical systems is refined into abstract orbit spaces, which is also a refinement of Reeb graphs of Hamiltonian flows with finitely many singular points on surfaces, and the CW decompositions which consist of the unstable manifolds of singular points for Morse flows on closed manifolds. 
Moreover, using the abstract orbit spaces, the author reconstructed some classes of flows from the time-one mappings. 

In this paper, interpreting concepts of recurrence, chain recurrence, and abstract orbits, we introduce abstract (weak) element spaces and Morse hyper-graphs of topological spaces and decompositions, show the existence of such hyper-graphs for topological spaces and invariant decompositions (e.g. compact foliated spaces), and reduce abstract (weak) element spaces into Morse hyper-graphs. 
Moreover, the abstract weak element space generalizes both the Reeb graph of a Morse function and the Morse graph for a flow, and refines an abstract cell complex.

The present paper consists of eight sections.
In the next section, recall some concepts of dynamical system to interpret such concepts in topological spaces and decompositions. 
In \S 3, we recall notions of combinatorics and topology, and the concepts of dynamical system are interpreted in topological spaces. In addition, the existence of Morse hyper-graphs is shown.  
In \S 4, we demonstrate that abstract element spaces can be reduced into Morse hyper-graphs. 
In \S 5, the correspondence with abstract element spaces and abstract cell complexes is constructed. More precisely, the abstract element space of a cell complex is an abstract cell complex with the specialization pre-order and its height. 
In \S 6, the recurrence of a flow on a locally compact Hausdorff space coincides with the recurrence for the orbit space. 
In \S 7, to apply results as above to decompositions (and especially compact foliated spaces), we also introduce concepts for decompositions on topological spaces as above, and show the existence of Morse hyper-graphs and the reducibility of abstract (weak) element spaces into such hyper-graphs. 
In addition, we demonstrate that the Reeb graph of a Morse function on a closed manifold is the abstract weak element space of the set of connected components of level sets as abstract multi-graphs. 
In the final section, some examples are illustrated to describe recurrent properties and the necessity of conditions. 

\section{Preliminaries from dynamical system's point of view}

We recall some concepts of dynamical systems to interpret such concepts in topological spaces and decompositions later. 

\subsection{Notion of dynamical systems}

A {\bf flow} is a continuous $\R$-action on a topological space.
Let $v : \R \times X \to X$ be a flow on a topological space $X$.
For $t \in \R$, define $v_t : X \to X$ by $v_t := v(t, \cdot )$.
For a point $x$ of $S$, we denote by $O(x)$ the orbit $\{ v_t(x) \mid t \in \R \}$ of $x$. 
\begin{definition}
A subset is {\bf invariant} {\rm(}or {\bf saturated}{\rm)} if it is a union of orbits. 
\end{definition}
The {\bf saturation} of a subset is the union of orbits intersecting it.
A nonempty closed invariant subset is {\bf minimal} if it contains no proper nonempty closed invariant subsets. 

\subsubsection{Orbit spaces and orbit class spaces}
The {\bf orbit space} $T/v$ of an invariant subset $T$ of $X$ is a quotient space $T/\sim$ defined by $x \sim y$ if $O(x) = O(y)$. 
Notice that an orbit space $T/v$ is the set $\{ O(x) \mid x \in T \}$ as a set. 
The {\bf (orbit) class} $\hat{O}$ of an orbit $O$ is the union of orbits each of whose orbit closure corresponds to $\overline{O}$ (i.e. $\hat{O} := \{ y \in S \mid \overline{O(y)} = \overline{O} \} $), where $\overline{O}$ is the closure of $O$.
The {\bf orbit class space} $T/\hat{v}$ of an invariant subset $T$ of $X$ is a quotient space $T/\approx$ defined by $x \approx y$ if $\overline{O(x)} = \overline{O(y)}$. 
Moreover, the orbit class space $T/\hat{v}$ is the set $\{ \hat{O}(x) \mid x \in T \}$ with the quotient topology. 
Denote by $\tau_v$ (resp. $\tau_{\hat{v}}$) the topology of the orbit space $X/v$ (resp. orbit class space $X/\hat{v}$). 
Note that the orbit class space $X/{\hat{v}}$ is  the $T_0$-tification (see the definition in \S \ref{topology}) of the orbit space $X/v$, and that 
$T/v$ (resp. $T/\hat{v}$) is a subset of the orbit space $X/v$ (resp. orbit class space $X/\hat{v}$).

\subsubsection{Properties of points and orbits}
A point $x$ of $S$ is {\bf singular} if $x = v_t(x)$ for any $t \in \R$ and is {\bf periodic} if there is a positive number $T > 0$ such that $x = v_T(x)$ and $x \neq v_t(x)$ for any $t \in (0, T)$.
A point is {\bf closed} if it is singular or periodic.
Denote by $\mathop{\mathrm{Sing}}(v)$ (resp. $\mathop{\mathrm{Per}}(v)$, $\bm{\mathop{\mathrm{Cl}}(v)}$) the set of singular (resp. periodic, closed) points. 
The {\bf $\bm{\omega}$-limit} (resp. {\bf $\bm{\alpha}$-limit}) set of a point $x$ is $\omega(x) := \bigcap_{n\in \mathbb{R}}\overline{\{v_t(x) \mid t > n\}}$ (resp.  $\alpha(x) := \bigcap_{n\in \mathbb{R}}\overline{\{v_t(x) \mid t < n\}}$).

\begin{definition}
A point $x$ of $S$ is {\bf recurrent} if $x \in \omega(x) \cup \alpha(x)$. 
\end{definition}

Denote by $\bm{\mathcal{R}(v)}$ (resp. $\bm{\mathrm{R}(v)}$) the set of recurrent (reps. non-closed recurrent) points and by $\bm{\mathrm{P}(v)}$ the set of non-recurrent points.
An orbit is closed (resp. recurrent) if it contains a closed (resp. recurrent) point. 

\begin{definition}
An orbit is {\bf proper} if it is embedded. 
\end{definition}

A point is proper if its orbit is proper. 

\begin{definition}
The {\bf derived set} of a subset $A$ is the set difference $\overline{A} - A$, where $B - C$ is used instead of the set difference $B \setminus C$ when $B \subseteq C$.
\end{definition}

\subsubsection{Characterizations of properness and recurrence}
We have the following observation. 

\begin{lemma}\label{lem:01}
An orbit $O$ of a flow on a topological space is proper if and only if the derived set $\overline{O} - O$ is closed. 
\end{lemma}

\begin{proof}
 Let $v$ be a flow on a topological space $X$. 
Suppose that $O \subseteq X$ is proper. 
Since $O$ is embedded, for any point $x \in O$, there is an open \nbd $U_x$ of $x$ such that $\overline{O} \cap U_x = O \cap U_x$. 
Then the union $U := \bigcup_{x \in O} U_x$ is an open \nbd of $O$ such that $\overline{O} \cap U = O$. 
Therefore the derived set $\overline{O} - O = \overline{O} \setminus U$ is closed. 
Conversely, suppose that $\overline{O} - O$ is closed. 
The complement $X - (\overline{O} - O)$ is an open \nbd of $O$ such that $\overline{O} \cap U = O$. 
This means that $O$ is embedded. 
\end{proof}

We have the following observation. 
\begin{lemma}\label{lem:02}
The following statements hold for a flow $v$ on a paracompact manifold: 
\\
{\rm(1)} The set $\mathrm{P}(v)$ of non-recurrent points is the set of non-closed proper points. 
\\
{\rm(2)} The set $\mathrm{R}(v)$ of non-closed recurrent points is the set of non-proper points. 
\\
{\rm(3)} The set $\mathcal{R}(v)$ of recurrent points is the set of closed or non-proper points. 
\end{lemma}

\begin{proof}
Notice that a closed orbit is proper and recurrent, and that $M = \Cv \sqcup \mathrm{P}(v) \sqcup \mathrm{R}(v)$, where $\sqcup$ denotes a disjoint union.  
We claim that any non-closed point is proper if and only if it is non-recurrent. 
Indeed, fix a non-closed point $x \in M$. 
By \cite[Corollary~3.4]{yokoyama2019properness}, the point $x$ is proper if and only if $O(x) = \hat{O}(x)$. 
From invariance of orbit and by definition of recurrence, if $x$ is non-recurrent then $O(x) = \hat{O}(x)$ and so proper.  
By \cite[Theorem~VI]{cherry1937topological}, the closure of a non-closed recurrent orbit $O$ of $v$ contains uncountably many non-closed recurrent orbits whose closures are $\overline{O}$. 
Then if $x$ is proper then $x$ is not recurrent.  
Therefore, the point $x$ is proper if and only if $x$ is non-recurrent. 

Since non-recurrent points are non-closed, the set $\mathrm{P}(v)$ of non-recurrent points is the set of non-closed proper points.  
Moreover, any non-closed point is non-proper if and only if it is recurrent. 
Since non-proper points are non-closed, the set of non-proper points is the set $\mathrm{R}(v)$ of non-closed recurrent points.  
\end{proof}


\subsubsection{Abstract weak orbits and abstract orbits}

Define an invariant subset $[x]$ for a flow on a topological space $X$ as follows:
\[
  [x] = \begin{cases}
    \text{the connected component of } \Sv \text{ containing } x & \text{if } x \in \mathop{\mathrm{Sing}}(v)  \\
    \text{the connected component of } \Pv \text{ containing } x & \text{if } x \in \mathop{\mathrm{Per}}(v)  \\
    \text{the connected component of } \\
     \{ y \in \mathrm{P}(v) \mid \alpha(x) = \alpha(y), \omega(x) = \omega(y) \}
     \text{ containing } x& \text{if } x \in \mathrm{P}(v) \\
    \{ y \in \mathrm{R}(v) \mid \alpha(x) = \alpha(y), \omega(x) = \omega(y) \} & \text{if } x \in \mathrm{R}(v)
  \end{cases}
\]
We call that $[x]$ is the {\bf abstract weak orbit} of $x$ (see \cite{yokoyama2021refinements} for details of properties of abstract orbits).
%
Similarly, define an invariant subset $\langle x \rangle$ of $X$ as follows:
\[
  \langle x \rangle = \begin{cases}
    \text{the connected component of } \Sv \text{ containing } x & \text{if } x \in \mathop{\mathrm{Sing}}(v)  \\
    \text{the connected component of } \Pv \text{ containing } x & \text{if } x \in \mathop{\mathrm{Per}}(v)  \\
    \text{the connected component of } \\
     \{ y \in \mathrm{P}(v) \mid \alpha(x) = \alpha(y), \omega(x) = \omega(y) \}
     \text{ containing } x& \text{if } x \in \mathrm{P}(v) \\
     \{ y \in \mathrm{R}(v) \mid \overline{O(x)} = \overline{O(y)} \} & \text{if } x \in \mathrm{R}(v)
  \end{cases}
\]
We call that $\langle x \rangle$ is the {\bf abstract orbit} of $x$ (see \cite{yokoyama2021refinements} for details of properties of abstract weak orbits).

\subsubsection{Chain recurrence}
Let $w \colon  \R \times M \to M$ be a flow on a metric space $(M,d)$.
For any $\varepsilon > 0$ and $T>0$, a pair $\{ (x_i)_{i=0}^{k+1}, (t_i)_{i=0}^{k} \}$ is an {\bf $\bm{(\varepsilon, T)}$-chain} from a point $x \in M$ to a point $y \in M$ if $x_0 = x$, $x_{p+1} = y$, $t_i > T$ and $d(w_{t_i}(x_i),x_{i+1}) < \varepsilon$ for any $i = 0, \ldots , k$.
Define a binary relation $\sim_{\mathop{CR}}$ on $M$ by $x \sim_{\mathop{CR}} y$ if for any $\varepsilon > 0$ and $T>0$ there is an $(\varepsilon, T)$-chain from $x$ to $y$.
\begin{definition}
A point $x \in M$ is {\bf chain recurrent} \cite{conley1978isolated} if $x \sim_{\mathop{CR}} x$.
\end{definition}

Denote by $\mathop{CR} (w)$ the set of chain recurrent points, called the {\bf chain recurrent set}.
It is known that the chain recurrent set $\mathop{CR} (w)$ of a flow on a compact metric space is closed and invariant and contains the non-wandering set $\Omega (w)$ \cite[Theorem~3.3B]{Conley1988}, and that connected components of $\mathop{CR} (w)$ are equivalence classes of the relation $\approx_{\mathop{CR}}$ on $\mathop{CR} (w)$ \cite[Theorem~3.3C]{Conley1988}, where $x \approx_{\mathop{CR}} y$ if $x \sim_{\mathop{CR}} y$ and $y \sim_{\mathop{CR}} x$.

\subsubsection{Morse graphs}
%
For a flow $v$ on a compact metric space $M$ with a set $\mathcal{M} = \{ M_i \}_{\lambda \in \Lambda}$ of disjoint compact invariant subsets, a finite directed graph $(V, D)$ with the vertex set $V := \{ M_i \mid 1 \leq i \leq n \}$, and with the {\bf directed edge set} $D := \{ (M_j, M_k) \mid D_{j,k} \neq \emptyset \}$ is a Morse graph of $\mathcal{M}$ if $M - \bigsqcup_i M_i = \bigsqcup D_{j,k}$, where 
\[
D_{j,k} := \left\{ x \in M - \bigsqcup_i M_i \mid \alpha(x) \subseteq M_j, \omega(x) \subseteq M_k \right\}
\]
and $\bigsqcup$ denotes a disjoint union. 
Then such a graph is denoted by $G_{\mathcal{M}}$ and called by the {\bf Morse graph} of $\mathcal{M}$. 

\begin{definition}
The graph $G_{\mathcal{M}}$ is the {\bf Morse graph} of $v$ if $\mathcal{M}$ is the set of connected components of the chain recurrent set $\mathrm{CR}(v)$. 
\end{definition}

Then denoted by $G_v$ the Morse graph of $v$. 
Similarly, if $\mathcal{M}$ is the set of connected components of the recurrent set $\mathcal{R}(v)$, then the graph $G_{\mathcal{M}}$ is denoted by $G_{\mathcal{R}(v)}$. 
Then a Morse graph $\mathcal{M}_{\mathcal{M}}$ of a flow $v$ is a quotient space of the orbit space with a directed structure. 
From \cite[Theorem 4.5]{yokoyama2021refinements}, any connected components of the chain recurrent point set are unions of abstract orbits.

\section{Preliminaries for the description of topological and combinatorial concepts}

In this section, we introduce topological invariants for topological spaces. 
To construct such invariants, we define some concepts (e.g. Morse graph, recurrence, ``chain recurrence'', ``abstract orbit'') for topological spaces from a dynamical system's point of view. 

\subsection{Fundamental notions of combinatorics and topology}

\subsubsection{Notion of combinatorics}

An {\bf abstract multi-graph} is a triple of sets $V, E$ and a mapping $A \colon E \to \{ \{ x, y \} \mid x,y \in V\}$. 
By a {\bf graph}, we mean a cell complex whose dimension is at most one and which is a geometric realization of an abstract multi-graph.
A {\bf (abstract) hyper-multi-graph} is a triple of sets $V, E$ and a mapping $A \colon E \to V^*$, where $V^*$ is the family of non-empty finite subsets of $V$. 
A {\bf (abstract) hyper-graph} is a pair of a set $V$ and a family $H \subseteq V^*$.

By a {\bf decomposition}, we mean a family $\mathcal{F}$ of pairwise disjoint nonempty subsets of a set $X$ such that $X = \bigsqcup \mathcal{F}$.
Since connectivity is not required, the sets of orbits of homeomorphisms are also decompositions.

A binary relation $\leq$ on a set $X$ is a {\bf pre-order} if it is reflexive (i.e. $a \leq a$ for any $a \in X$) and transitive (i.e. $a \leq c$ for any $a, b, c \in X$ with $a \leq b$ and $b \leq c$).
For a pre-order $\leq$, the inequality $a<b$ means both $a \leq b$ and $a \neq b$.
A pre-order order $\leq$ is a {\bf total order} (or linear order) if either $a < b$ or $b < a$ for any  points $a \neq b$.
A {\bf chain} is a totally ordered subset of a pre-ordered set with respect to the induced order.
Let $(X, \leq)$ be a pre-ordered set.
Define the height $\mathop{\mathrm{ht}} (x)$ of $x \in X$ by $\mathop{\mathrm{ht}} (x) := \sup \{ |C| - 1 \mid C :\text{chain containing }x \text{ as the maximal point}\}$. 
Define the height of the empty set as $-1$. 
The {\bf height} $\mathop{\mathrm{ht}} (A)$ of a nonempty subset $A \subseteq X$ is defined by $\mathop{\mathrm{ht}} (A) := \sup_{x \in A} \mathop{\mathrm{ht}} (x)$.
A subset $A \subset X$ is a downset if $b \in A$ for any point $a \in A$ and any point $b \in X$ with $b \leq a$.

\subsubsection{Notion of topology}\label{topology}

Let $(X, \tau)$ be a topological space. 
A point $x \in X$ is $\bm{T_0}$ (or {\bf Kolmogorov}) if for any point $y \in X - \{ x \}$, there is an open subset $U$ of such that $\{x, y \} \cap U$ is a singleton,  $\bm{T_1}$ (or {\bf closed}) if the singleton $\{ x \}$ is closed, and $\bm{T_D}$ \cite{aull1962separation} if the derived set $\overline{\{ x \}} - \{ x \}$ is a closed subset. 
Here a singleton is a set consisting of a point.
A topological space is $T_0$ (resp. $T_1$, $T_D$, etc.) if each point is $T_0$ (resp. $T_1$, $T_D$, etc.).

\begin{definition}
The {\bf class $\bm{\hat{x}}$} of a point $x \in X$ is defined by a subset $\{ y \in X \mid \overline{\{ x \}} = \overline{\{ y \}}  \}$.
\end{definition} 
Then the set $\hat{X} := \{ \hat{x} \mid x \in X \}$ of classes is a decomposition of $X$ and is a $T_0$ space as a quotient space, which is called the {\bf $\bm{T_0}$-tification} (or {\bf Kolmogorov quotient}) of $X$.
\begin{definition}
A subset is {\bf invariant} if it is a union of classes. 
\end{definition}

Notice that any closed subset is invariant. 
%
%
The {\bf specialization pre-order} $\leq_{\tau}$ on a topological space $(X, \tau)$ is defined as follows: $ x \leq_{\tau} y $ if  $ x \in \overline{\{ y \}}$. 
The heights of a point, a subset, and a topological space are defined as the heights with respect to the specialization pre-order. 
For any $k \in \Z_{\geq 0}$ and a topological space $X$, denote by $X_k$ the set of height $k$ points of $X$ and by $X_{\leq k}$ the set of points of $X$ whose heights are less than or equal to $k$.

\subsection{Topological concepts for topological spaces from Dynamical systems}

Let $(X, \tau)$ be a topological space. 
By Lemma~\ref{lem:01}, we define properness as follows. 
\begin{definition}
A point of $X$ is {\bf proper} if it is $T_D$. 
\end{definition}
Denote by $\bm{\Ct}$ (resp. $\bm{\mathrm{P}(\tau)}$, $\bm{\mathrm{R}(\tau)}$) the set of closed points (resp. non-closed proper points, non-proper points). 
By definition, we have $X = \Ct \sqcup \mathrm{P}(\tau) \sqcup \mathrm{R}(\tau)$. 
From Lemma~\ref{lem:02}, we call that a point $x \in X$ is {\bf ($\bm{\tau}$-)recurrent} if it is either $T_0$ or non-$T_D$. 
Denote by $\bm{\mathcal{R}(\tau)}$ the set of recurrent points. 
Then $\mathcal{R}(\tau) = \Ct \sqcup \mathrm{R}(\tau)$ and $X = \mathrm{P}(\tau) \sqcup \mathcal{R}(\tau)$. 
We will show that an orbit of a flow $v$ on a locally compact Hausdorff space $X$ is recurrent if and only if it is $\tau_v$-recurrent, where $\tau_v$ is the quotient topology of the orbit space $X/v$ (see Theorem~\ref{lem:23}). 
\begin{definition}
Define an {\bf abstract element} $\bm{\langle x \rangle}$ for an element $x \in X$:   
\[
  \langle x \rangle := \begin{cases}
    \text{the connected component of} \\
     \{ x' \in X \mid \overline{\{ x \}} - \{ x \} = \overline{\{ x' \}} - \{ x' \} \} 
     \text{ containing } x & \text{if } x \in \mathrm{Cl}(\tau) \sqcup \mathrm{P}(\tau) \\
         \text{the connected component of} \\
   \{ x' \in X \mid \overline{\{ x \}} = \overline{\{ x' \}} \}  \text{ containing } x  & \text{if } x \in \mathrm{R}(\tau) 
  \end{cases}
\]
\end{definition}

We have the following observation. 

\begin{lemma}\label{lem:ch_ab_ele}
The abstract elements form a decomposition and satisfy the following property: 
\[
  \langle x \rangle := \begin{cases}
    \text{the connected component of } \mathrm{Cl}(\tau) \text{ containing } x & \text{if } x \in \mathrm{Cl}(\tau) \\
    \text{the connected component of} \\
     \{ x' \in \mathrm{P}(\tau) \mid \overline{\{ x \}} - \{ x \} = \overline{\{ x' \}} - \{ x' \} \} 
     \text{ containing } x & \text{if } x \in \mathrm{P}(\tau) \\
         \text{the connected component of} \\
  \{ x' \in \mathrm{R}(\tau) \mid \overline{\{ x \}} = \overline{\{ x' \}} \}  \text{ containing } x  & \text{if } x \in \mathrm{R}(\tau) 
  \end{cases}
\]
\end{lemma}

\begin{proof}
Let $(X, \tau)$ be a topological space. 
Fix a point $x \in X$. 
If $x \in \mathrm{Cl}(\tau)$, then $\langle x \rangle$ is the connected component of $\{ x' \in X \mid \overline{\{ x \}} - \{ x \} = \overline{\{ x' \}} - \{ x' \} = \emptyset \} = \Ct$. 

Suppose that $x \in \mathrm{R}(\tau)$. 
If there is a point $y \in \langle x \rangle \cap \Ct$, then $x \in \overline{\{ y \}} = \{ y \} \subseteq \Ct$, which contradicts $x \in \mathrm{R}(\tau) = X - (\Ct \sqcup \mathrm{P}(\tau))$. 
Thus $\langle x \rangle \cap \Ct = \emptyset$. 
We claim that $\langle x \rangle \cap \mathrm{P}(\tau) = \emptyset$. 
Indeed, assume that there is a point $y \in \langle x \rangle \cap \mathrm{P}(\tau)$. 
Since $x \in \mathrm{R}(\tau)$, we have $x \in \overline{\{ y \}} - \{ y \}$. 
Since $\overline{\{ y \}} - \{ y \}$ is closed, we have $\overline{\{ y \}} = \overline{\{ x \}} \subseteq \overline{\{ y \}} - \{ y \}$, which is a contradiction. 
Thus $\langle x \rangle$ is the connected component of $\{ x' \in \mathrm{R}(\tau) \mid \overline{\{ x \}}  = \overline{\{ x' \}} \}$. 

Suppose that $x \in \mathrm{P}(\tau)$. 
Then the derived set $\overline{\{ x \}} - \{ x \} \neq \emptyset$ is closed, and the abstract element $\langle x \rangle$ is the connected component of $\{ x' \in X \mid \overline{\{ x \}} - \{ x \} = \overline{\{ x' \}} - \{ x' \} \neq \emptyset \}$. 
Therefore $\langle x \rangle \cap \Ct = \emptyset$. 
We claim that $\langle x \rangle \cap \mathrm{R}(\tau) = \emptyset$. 
Indeed, assume that there is a point $y \in \langle x \rangle \cap \mathrm{R}(\tau)$. 
Then $y \in \overline{\{ x \}} - \{ x \}$. 
Since $\overline{\{ x \}} - \{ x \}$ is closed, the set difference $\overline{\{ y \}} - \{ y \} = \overline{\{ x \}} - \{ x \}$ is closed and so $y \in \mathrm{P}(\F)$, which contradicts $y \in \mathrm{R}(\tau)$. 
Thus $\langle x \rangle$ is the connected component of $\{ x' \in \mathrm{P}(\tau) \mid \overline{\{ x \}} - \{ x \} = \overline{\{ x' \}} - \{ x' \} \}$. 
\end{proof}

%
An abstract element is closed (resp. recurrent, proper, etc.) if it is contained in $\Ct$ (resp. $\mathcal{R}$, $\Ct \sqcup \mathrm{P}(\tau)$, etc.). 
Define the {\bf abstract element space} $X/\langle \tau \rangle$ as a quotient space $X/\sim_{\langle \tau \rangle}$ defined by $x \sim_{\langle \tau \rangle} y$ if $\langle x \rangle = \langle y \rangle$.

\subsubsection{Quasi-recurrence of points}
We will define an analogous concept of chain recurrence. 
However, the absence of metrics obstructs a direct interpretation of chain recurrence. 
To construct ``chain recurrence'', recall facts that any points in the orbit closure for a flow are chain recurrent, and that any connected components of the chain recurrent point set are unions of abstract orbits. 
To achieve the analogous concept, we define non-maximal points as follows. 
\begin{definition}
A point $x \in X$ is {\bf non-maximal} if there is a point $y \in X$ such that $\overline{\{ x \}} \subsetneq \overline{\{ y \}}$.
\end{definition} 
Note that a point is non-maximal if and only if it is not maximal with respect to the specialization pre-order of $\tau$. 
A point is maximal if it is not non-maximal. 
Denote by $\max \tau$ the set of maximal points. 
%
Then we define ($\tau$-)quasi-recurrence for a topological space, which is the analogous concept of chain recurrence, as follows. 
\begin{definition}
A point $x \in X$ is {\bf {\rm(}$\bm{\tau}$-{\rm)}quasi-recurrent} if $\langle x \rangle$ contains a recurrent or non-maximal point. 
\end{definition}

Denote by $\bm{\mathcal{Q}(\tau)}$ the set of quasi-recurrent points. 
Notice that $\mathcal{Q}(\tau)$ is a union of abstract elements but is not closed in general (see the example in \S~\ref{sec:01}). 

\subsubsection{Morse hyper-graph for a topological space}

We define the Morse hyper-graph $(V, H)$ for a topological space as follows: 
Let $(X, \tau)$ be a topological space with a set $\mathcal{X} = \{ X_{\lambda} \}_{\lambda \in \Lambda}$ of disjoint invariant non-empty subsets $X_\lambda$ ($\lambda \in \Lambda$). 
For any $I \subseteq \lambda$, define a {\bf hyper-edge} $H_{I}$ as follows: For any point $x \in X - \bigsqcup_{\lambda \in \Lambda} X_\lambda$, we say that $x \in H_{I}$ if there are disjoint invariant non-empty subsets $(C_i)_{i \in I}$ of the subset $\overline{\{ x \}} - \{ x \}$ such that $C_i \subseteq X_i$ and $\overline{\{ x \}} - \{ x \} = \bigsqcup_{i \in I} C_i$. 
Put $V := \{ X_{\lambda} \mid \lambda \in \Lambda \}$ and $H := \{ \{ X_i \}_{i \in I} \mid H_{I} \neq \emptyset, I \subseteq \Lambda \}$. 
\begin{definition}
A hyper-graph $\mathcal{G}_{\mathcal{X}} := (V, H)$ is a {\bf Morse hyper-graph} of $\mathcal{X}$ if $X - \bigsqcup_{\lambda \in \Lambda} X_\lambda = \bigsqcup_{I \subseteq \Lambda} H_{I}$. 
\end{definition}
\begin{definition}
The hyper-graph $\mathcal{G}_{\mathcal{X}}$ is the {\bf Morse hyper-graph} of $\tau$ if $\mathcal{X}$ is the set of connected components of the set $\mathcal{Q}(\tau)$ of quasi-recurrent points. 
\end{definition}


Notice that we can similarly define the Morse hyper-multi-graph. 
We have the following statement.

\begin{lemma}\label{lem:ch_fol_rec_chain_top}
The following statements hold for a topological space $(X, \tau)$: 
\\
{\rm(1)} $\mathcal{R} (\tau) \sqcup (\bigcup_{x \in X} (\overline{\{ x \}} - \{ x \}) \cap \mathrm{P}(\tau)) \subseteq \mathcal{Q}(\tau)$. 
\\
{\rm(2)} $X - \mathcal{Q}(\tau) \subseteq \mathrm{P}(\tau) \cap \max \tau$. 
\\
{\rm(3)} The set $\mathcal{Q}(\tau)$ of quasi-recurrent points is invariant. 
\\
{\rm(4)} For any element $ x \in \mathrm{P}(\tau) \cap \max \tau$, each connected component of $\overline{\{ x \}} - \{ x \}$ is a closed subset contained in a connected component of $\mathcal{Q}(\tau)$. 
\end{lemma}

\begin{proof}
The non-minimal property implies that the set of non-maximal points is invariant. 
By definition of quasi-recurrence, the set $\mathcal{Q}(\tau)$ of quasi-recurrent points is invariant such that $\bigcup_{x \in X} \overline{\{ x \}} - \{ x \} \subseteq \mathcal{Q}(\tau)$, $X - \mathcal{Q}(\tau) \subseteq \max \tau$, and $X - \mathrm{P}(\tau) = \mathcal{R} (\tau) \subseteq \mathcal{Q}(\tau)$. 
Then $\mathcal{R} (\tau) \sqcup (\bigcup_{x \in X} (\overline{\{ x \}} - \{ x \}) \cap \mathrm{P}(\tau)) \subseteq \mathcal{Q}(\tau)$ and $X - \mathcal{Q}(\tau) \subseteq \mathrm{P}(\tau) \cap \max \tau$.
 
For any point $x \in \mathrm{P}(\tau) \cap \max \tau$, 
the derived subset $\overline{\{ x \}} - \{ x \}$ is a closed subset contained in $X - \max \tau \subseteq \mathcal{Q}(\tau)$ and so each connected component of $\overline{\{ x \}} - \{ x \}$ is contained in a connected component of $\mathcal{Q}(\tau)$.  
\end{proof}

The equality in Lemma~\ref{lem:ch_fol_rec_chain_top} {\rm(1)} does not hold in general (see Example in \S \ref{ex:04}). 
Lemma~\ref{lem:ch_fol_rec_chain_top} implies the following existence of Morse hyper-graph for a topological space. 

\begin{lemma}\label{lem:ex_MH_top_02}
Let be $(X, \tau)$ a topological space, $\mathcal{Q}' \subseteq X$ an invariant subset containing $\mathcal{Q}(\tau)$, and $\mathcal{X} = \{ X_{\lambda} \}_{\lambda \in \Lambda}$ the set of connected components of $\mathcal{Q}'$. 
Then the Morse hyper-graph $\mathcal{G}_\mathcal{X}$ exists. 
\end{lemma}

\begin{proof}
Since any connected component of a subset is closed in the subset, any connected component $X_\lambda$ is invariant. 
Therefore $\mathcal{X}$ is a set of disjoint invariant subsets. 
Lemma~\ref{lem:ch_fol_rec_chain_top} implies that $X - \bigsqcup_{\lambda \in \Lambda} X_\lambda = X - \mathcal{Q}' \subseteq X - \mathcal{Q}(\tau) \subseteq \mathrm{P}(\tau) \cap \max \tau$. 
Fix any point $ x \in X - \bigsqcup_{\lambda \in \Lambda} X_\lambda \subseteq X - \mathcal{Q}(\tau) \subseteq \mathrm{P}(\tau) \cap \max \tau$. 
By Lemma~\ref{lem:ch_fol_rec_chain_top}, the derived subset $\overline{\{ x \}} - \{ x \} \subseteq \mathcal{Q}(\tau) = \bigsqcup_{\lambda \in \Lambda} X_\lambda$ is closed and so the intersection $C_\lambda := (\overline{\{ x \}} - \{ x \}) \cap X_\lambda$ is invariant.
Set $I := \{ \lambda \in \Lambda \mid C_\lambda \neq \emptyset \}$. 
Then $x \in  H_I$. 
This means that $X - \bigcup_{\lambda \in \Lambda} X_\lambda = \bigsqcup_{I \subseteq \Lambda} H_{I}$. 
\end{proof}

The previous lemmas imply the following statement. 

\begin{corollary}\label{lem:ex_MH_top}
The Morse hyper-graph for a topological space exists.
\end{corollary}

\section{Reductions to Morse hyper-graphs for topological spaces}


The Morse hyper-graph of a topological space can be obtained as a quotient space of the abstract element space as follows.

\begin{theorem}\label{th:Morse_reduction}
The Morse hyper-graph $\mathcal{G}_{\tau}$ for a topological space $(X, \tau)$ is a quotient space of the abstract element space $X/\langle \tau \rangle$.
\end{theorem}

\begin{proof}
Let $\mathcal{G}_{\tau} = (V, H)$ be the Morse hyper-graph for a topology $\tau$ with $V = \{ X_\lambda \}_{\lambda \in \Lambda}$ and $H = \{ \{ X_i \}_{i \in I} \mid H_{I} \neq \emptyset, I \subseteq \Lambda \}$.
Then $\mathcal{X} = \{ X_\lambda \}_{\lambda \in \Lambda}$ is the set of connected components $X_\lambda$ of $\mathcal{Q}(\tau)$. 
By $\mathcal{R} (\tau) \subseteq \mathcal{Q}(\tau) = \bigsqcup_{\lambda \in \Lambda} X_\lambda$, since any abstract elements are connected, any recurrent abstract elements are contained in some $X_\lambda$. 

Fix any point $x \in \mathrm{P}(\tau)$. 
Then the abstract element $\langle x \rangle$ is a connected component of $\{ x' \in \mathrm{P}(\tau) \mid \overline{\{ x\}} - \{ x \} = \overline{\{ x' \}} - \{x' \} \neq \emptyset \}$ and the derived set $\overline{\{ x\}} - \{ x \}$ is closed. 
Suppose that $x$ is quasi-recurrent. 
Then $\langle x \rangle$ is connected and contained in $\mathcal{Q}(\tau)$. 
Therefore there is an index $i \in \Lambda$ such that $\langle x \rangle \subseteq X_i$. 
Suppose that $x$ is not quasi-recurrent. 
Then $\langle x \rangle \subseteq X - \mathcal{Q}(\tau) = X - \bigsqcup_{\lambda \in \Lambda} X_\lambda  = \bigsqcup_{I \subseteq \Lambda} H_{I}$. 
Fix any $x' \in \langle x \rangle$. 
Then $\overline{\{ x\}} - \{ x \} = \overline{\{ x' \}} - \{x' \} \subseteq \mathcal{Q}(\tau) = \bigsqcup_{\lambda \in \Lambda} X_\lambda$. 
Put $C_i := X_i \cap (\overline{\{ x\}} - \{ x \}) = X_i \cap (\overline{\{ x' \}} - \{ x' \})$. 
Define $I := \{ i \in \Lambda \mid C_i \neq \emptyset \}$. 
Since any connected components of a subset are closed in the subset, the intersections $C_i$ for any $i \in I$ are invariant such that $C_i \subseteq X_i$ and $\overline{\{ x' \}} - \{ x' \} = \overline{\{ x \}} - \{ x \} =  (\overline{\{ x \}} - \{ x \}) \cap \bigsqcup_{\lambda \in \Lambda} X_\lambda = \bigsqcup_{i \in I} C_i$. 
Then $x', x \in H_I$ and so $ \langle x \rangle \subseteq H_I$. 
\end{proof}

In general, the Morse hyper-graph $\mathcal{G}_{\mathcal{X}}$ even for $\mathcal{X}$ as in Lemma~\ref{lem:ex_MH_top_02} is not a quotient space of the abstract element space (see the example in \ref{ex:mhg}).

\section{Generalization of abstract cell complex structures}

This section shows that the abstract element space is a generalization of an abstract cell complex structures. 

\subsection{Generalization of abstract cell complex structure}

We recall concepts of cell complexes. 


\begin{definition}
A decomposition $\{ e_\lambda \}_{\lambda \in \Lambda}$ on a Hausdorff space $X$ is a {\bf cell complex} if it satisfies the following two conditions: 
\\
{\rm(1)} Any element $e_\lambda$ is a $k$-dimensional ball for some $k \in \Z_{\geq 0}$. In this case, the element $e_\lambda$ is called a $k$-cell, and the integer $k$ is called the dimension of $e_\lambda$. 
\\
{\rm(2)} For any $k$-cell $e_\lambda$, we have $\overline{e_\lambda} - e_\lambda \subseteq X^{\leq k-1}$, where $X^{\leq k-1}$ is the union of cells whose dimension is less than $k$. 
\end{definition}

\begin{definition}
For a set $S$ with a transitive relation $\prec$ and a function $\dim \colon S \to \Z_{\geq 0}$, the triple $(S, \prec, \dim )$ is an {\bf abstract cell complex} if $x \prec y$ implies $\dim x < \dim y$. 
\end{definition}
Then $\dim x$ is called the {\bf dimension} of $x$, and $x$ is called a {\bf cell}. 
A {\bf $\bm{k}$-cell} is a cell whose dimension is $k$. 
%

\subsubsection{Correspondence with abstract element spaces and abstract cell complexes}
We have the following correspondence. 

\begin{proposition}\label{prop:abst_cell}
The abstract element space of a cell complex is an abstract cell complex with the specialization pre-order and its height. 
\end{proposition}

\begin{proof}
Let $\F$ be a cell complex on a topological space $X$. 
By definition of cell complex, the abstract element space $X/\langle \F \rangle$ is $\F$ as a set, and the height $\mathrm{ht}(L)$ of any element $L$ of the quotient space $X/\langle \F \rangle$ with respect to the specialization pre-order $\leq$ corresponds to the dimension of the cell of $\F$. 
Then the tuple $(S, \prec, \dim ) := (X/\langle \F \rangle, \leq, \mathrm{ht})$ is an abstract cell complex. 
\end{proof}

\section{Correspondence of recurrences}

In this section, we demonstrate the correspondence of recurrences for flows and topological spaces.
By definitions,  the $\alpha$-limit set and the $\omega$-limit set of $x$ are closed and invariant.
We have the following observation. 
\begin{lemma}\label{lem:decomp_limit}
Let $v$ be a flow on a Hausdorff space $X$ and $x \in X$ a point. 
Then $\overline{O(x)} = \alpha(x) \cup O(x) \cup \omega(x)$.
\end{lemma}

\begin{proof}
By definition of $\alpha$-limit set and $\omega$-limit set, the orbit closure $\overline{O(x)}$ contains $\alpha(x) \cup O(x) \cup \omega(x)$.
Define a continuous mapping $v_x  \colon \R \to X$ by $v_x :=v( \cdot, x)$.
Since a closed interval $I \subset \R$ is compact, the image $v_x(I)$ is compact.
The Hausdorff separation axiom of $X$ implies that the image of a closed interval by $v_x$ is closed.

Assume that there is a point $y \in \overline{O(x)} - (\alpha(x) \cup O(x) \cup \omega(x))$.
By $y \notin \alpha(x) \cup \omega(x)$, there is a number $T > 0$, such that $y \notin \overline{v_x(\R_{<-T})}$ and $y \notin \overline{v_x(\R_{>T})}$.
Then $y \notin \overline{v_x(\R - [-T,T])}$ and so there is an open \nbd $U$ of $y$ such that $U \cap v_x(\R - [-T,T]) = \emptyset$.
Since the image $v_x([-T,T]) \subseteq O(x)$ of the interval $[-T,T]$ is closed and $y \notin O(x)$, the difference $V := U \setminus v_x([-T,T])$ is an open \nbd of $y$.
Then $V \cap O(x) = (U \setminus v_x([-T,T])) \cap O(x) = U \cap (O(x) - v_x([-T,T])) = U \cap v_x(\R - [-T,T]) = \emptyset$ and so $y \notin \overline{O(x)}$, which contradicts $y \in \overline{O(x)}$.
Thus $\overline{O(x)} = \alpha(x) \cup O(x) \cup \omega(x)$.
\end{proof}

We have the following relation between recurrences. 

\begin{theorem}\label{lem:23}
Let $v$ be a flow on a Hausdorff space $X$ and $\tau_{X/v}$ the quotient topology of the orbit space $X/v$. 
A $\tau_{X/v}$-recurrent orbit is $v$-recurrent. 
If $X$ is locally compact, then the converse holds. 
\end{theorem}

\begin{proof}
Suppose that $O(x)$ is an $\tau_{X/v}$-recurrent orbit. 
Then either $O(x)$ is closed or $\overline{O(x)} - O(x)$ is not closed. 
If $O(x)$ is closed, then $x \in \alpha(x) \cup \omega(x)$ trivially. 
Thus we may assume that $\overline{O(x)} - O(x)$ is not closed. 
By Lemma~\ref{lem:decomp_limit}, this implies that $x \in \overline{\overline{O(x)} - O(x)} \subseteq \overline{\alpha(x) \cup \omega(x)} = \alpha(x) \cup \omega(x)$. 

Suppose that $X$ is locally compact and $O(x)$ is a $v$-recurrent orbit. 
We may assume that $O(x)$ is not closed. 
We claim that $\overline{O(x)} - O(x)$ is not closed. 
Indeed, assume $\overline{O(x)} - O(x)$ is closed. 
Applying the Baire category theorem, 
since $\overline{O(x)}$ is locally compact Hausdorff and so Baire, and since each open subset of a Baire space is a Baire space, the orbit $O(x)$ is Baire. 
Let $U_n := v(\R - [n, n+1], x) \subset O(x)$ for $n \in \Z$. 
Because $O(x)$ is $v$-recurrent, each $U_n$ is open dense in $O(x)$. 
Since $O(x)$ is Baire, we have $\bigcap_n U_n$ is dense, which contradicts to the definition of $U_n$. 
By $\overline{O(x)} - O(x) \subseteq \alpha(x) \cup \omega(x)$, the non-closedness of $\overline{O(x)} - O(x)$ and the closed invariance of $\alpha(x) \cup \omega(x)$ imply that $O(x) \subset \alpha(x) \cup \omega(x)$. 
Therefore $O(x)$ is $\tau_{X/v}$-recurrent.  
\end{proof}

There is a transitive flow  on a metrizable space such that each non-singular point is not $\tau_{X/v}$-recurrent but $v$-recurrent (see an example in \S \ref{ex:02} for details).
 
%

\section{Application to decompositions and foliated spaces}

To apply previous results to decompositions and foliated spaces, we recall concepts of decompositions and introduce topological invariants for decompositions on topological spaces. 

\subsection{Notion of decompositions}
By a decomposition, we mean a family $\mathcal{F}$ of pairwise disjoint nonempty subsets of a set $X$ such that $X = \bigsqcup \mathcal{F}$, where $\bigsqcup$ denotes a disjoint union. 
Let $\F$ be a decomposition on a topological space $X$. 
%
%
For $x \in X$, denote by $\F(x)$ the element containing $x$. 
An element $L \in \F$ is {\bf closed} if $\overline{L}= L$. 
\begin{definition}
An element $L \in \F$ is {\bf proper} if the derived set $\overline{L}- L$ is closed. 
\end{definition}
Denote by $\mathop{\mathrm{Cl}}(\F)$ (resp. $\mathrm{P}(\mathcal{F})$, $\mathrm{R}(\mathcal{F})$) the union of closed (resp. non-closed proper, non-closed non-proper) elements. 
Set $\mathcal{R}(\mathcal{F}):= \mathop{\mathrm{Cl}}(\F) \sqcup \mathrm{R}(\mathcal{F}) = X - \mathrm{P}(\mathcal{F})$. 

\begin{definition}
A subset $A$ of $X$ is {\bf $\F$-saturated} {\rm(}or {\bf $\F$-invariant}{\rm)} if $\F(A) = A$, where $\F(A) = \bigcup_{x \in A} \F(x)$. 
\end{definition}

The subset $\F(A)$ is called the {\bf saturation} of $A$. 

\begin{definition}
$\F$ is said to be $T_{\sigma}$ (resp. recurrent, etc.) if the quotient space $X/{\F}$ of $\F$ is $T_{\sigma}$ (resp. recurrent, etc.). 
\end{definition}

Define the {\bf class element} $\hat {L}$ of an element $L \in \F$ by $\hat {L} := \{ x \in X \mid \overline{\F(x)} = \overline{L} \}$. 
Then the family $\hat {\F} := \{ \hat{L}  \mid L \in \F \}$ is a decomposition, called the {\bf class decomposition} of $\F$. 
The quotient space of the topological space $X$ by the decomposition $\F$ is called the {\bf class (element) space} of $\F$ and is denoted by $X/\hat{\F}$. 
%
%
%
%
Denote by $\tau_{\F} := \{ \F(U) \mid U \in \tau \}$ the set of $\F$-saturation of open subsets. 
\begin{definition}
A decomposition $\F$ is {\bf invariant} if $\tau_{\F} \subseteq \tau$. 
\end{definition}

We will show that if $\F$ is either a foliated space or a continuous action of a topological group then $\F$ is invariant (see Proposition~\ref{prop:02}). 
For an invariant decomposition $\F$, the set $\tau_{\F}$ of saturations of open subsets becomes a topology and called the {\bf saturated topology} on $X$.
In general, the set $\tau_{\F}$ of $\F$-invariant subsets is not a topology on the quotient space $X/\F$ even if $\overline{L}$ is $\F$-invariant for any $L \in \F$ (see an example in \S \ref{ex:01} for details).  

\subsubsection{Properties of decompositions}

We observe the following statements.

\begin{lemma}\label{lem:saturated}
The $\F$-invariant subsets $\mathop{\mathrm{Cl}}(\F)$, $\mathrm{P}(\mathcal{F})$, and $\mathrm{R}(\mathcal{F})$ of an invariant decomposition $\F$ on a topological space are $\hat{\F}$-invariant. 
\end{lemma}

\begin{proof}
Let  $\F$ be an invariant decomposition on a topological space $X$. 
For any point $x \in \mathop{\mathrm{Cl}}(\F)$, we have $\hat{\F}(x) = \{ y \in X \mid \overline{\F(y)} = \overline{\F(x)} = \F(x) \} \subseteq \mathop{\mathrm{Cl}}(\F)$.
This means that $\mathop{\mathrm{Cl}}(\F)$ is $\hat{\F}$-invariant. 
Fix a point $x \in \mathrm{P}(\mathcal{F})$. 
Since the derive set $\overline{\F(x)} - \F(x)$ is closed and $\F$-invariant, if there is a point $y \in \hat{\F}(x) - \F(x)$, then $\overline{\F(x)} = \overline{\F(y)} \subseteq \overline{\F(x)} - \F(x)$, which is a contradiction. 
Thus $\hat{\F}(x) = \F(x) \subseteq \mathrm{P}(\mathcal{F})$. 
Since $\mathop{\mathrm{Cl}}(\F) \sqcup \mathrm{R}(\mathcal{F}) = X - \mathrm{P}(\mathcal{F})$ is $\hat{\F}$-invariant, so is the complement $\mathrm{P}(\mathcal{F})$. 
\end{proof}

\begin{lemma}\label{lem:001}
The following statements are equivalent for a decomposition $\F$ on a topological space: 
\\
{\rm(1)} The decomposition $\F$ is invariant. 
\\
{\rm(2)} The closure of any $\F$-invariant subset is $\F$-invariant. 
\end{lemma}

%
%
%

\begin{proof} 
Let  $\F$ be a decomposition on a topological space $(X, \tau)$. 
Suppose that $\tau_{\F} \subseteq \tau$. 
Assume that there is an $\F$-invariant subset $A$ whose closure is not $\F$-invariant. 
Then there is a point $x \in \F(\overline{A}) - \overline{A}$. 
Put $U := X - \overline{A}$. 
This $U$ is an open \nbd of $x$. 
Since $\tau_{\F} \subseteq \tau$, the saturation $V := \F(U)$ is open with $V \cap A = \emptyset$. 
Then $V \cap \overline{A} = \emptyset$. 
Since $V$ is $\F$-invariant, we obtain $V \cap \F(\overline{A}) = \emptyset$. 
This contradicts that $V$ is a \nbd of $x \in \F(\overline{A})$. 
Thus $\overline{A}$ is $\F$-invariant for any $\F$-invariant subset $A \subseteq X$. 

Conversely, suppose that the closure of an $\F$-invariant subset is $\F$-invariant. 
Fix any open subset $B$. 
Set $F := X - \F(B)$. 
Since $B \cap F = \emptyset$, we have $B \cap  \overline{F} = \emptyset$. 
The hypothesis implies that $\overline{F}$ is $\F$-invariant. 
Then $\F(B) \cap  \overline{F} = \emptyset$ and so $\F(B) \cap  \overline{X - \F(B)} = \emptyset$. 
This implies that $X - \F(B)$ is closed and so $\F(B)$ is open. 
Therefore $\tau_{\F} \subseteq \tau$. 
\end{proof}

We have the following correspondence. 

\begin{lemma}\label{cor:correspondence}
Let $\F$ be an invariant decomposition on a topological space $X$ and $p \colon X \to X/\F$ the quotient map. 
The following statements hold: 
\\
{\rm(1)} $\{ \bigcup \mathcal{U} \mid \mathcal{U} \in \tau_{X/\F} \} =\{ p^{-1} (\mathcal{U}) \mid \mathcal{U} \in \tau_{X/\F} \} = \tau_{\F}$. 
\\
{\rm(2)} $\{ p(U) \mid U \in \tau_{\F} \} = \tau_{X/\F}$. 
\\
{\rm(3)} The quotient map $p$ induces the canonical bijection $\tau_{\F} \to \tau_{X/\F}$ by $U \mapsto p(U)$. 
\end{lemma}

\begin{proof}
The union of any elements of the quotient topology $\tau_{X/\F}$ is invariant open and so $\{ \bigcup \mathcal{U} \mid \mathcal{U} \in \tau_{X/\F} \} =\{ p^{-1} (\mathcal{U}) \mid \mathcal{U} \in \tau_{X/\F} \} \subseteq \tau_{\F}$.  
From $\mathcal{U} = p (p^{-1} (\mathcal{U}))$ for any $\mathcal{U} \in \tau_{X/\F}$, we have $\tau_{X/\F}  \subseteq \{ p (U) \mid U \in \tau_{\F} \}$. 
By $U = p^{-1} (p (U))$ for any $U \in \tau_{\F}$, since any elements of $\tau_{\F}$ are invariant open, we obtain $\{ p(U) \mid U \in \tau_{\F} \} \subseteq \tau_{X/\F}$.  
Thus $\{ p(U) \mid U \in \tau_{\F} \} = \tau_{X/\F}$. 
Because $p^{-1}(p(U)) = U$ for any $U \in \tau_{\F}$, we have $\tau_{\F} = \{ p^{-1} (\mathcal{U}) \mid \mathcal{U} \in \tau_{X/\F} \}$. 
This means that $p$ induces the canonical bijection $\tau_{\F} \to \tau_{X/\F}$. 
\end{proof}

From the previous lemma, we can identify the saturated topology of an invariant decomposition on a topological space as the topology on the decomposition space.
%

\subsection{Topological concepts for decompositions from dynamical systems}

We define topological concepts for decompositions from dynamical systems as the case for topological spaces. 

\subsubsection{Abstract weak element and abstract element for a decomposition}

Let $\mathcal{F}$ be a decomposition on a topological space $X$ and $q \colon X \to X/\F$ the quotient map. 
\begin{definition}
Define an {\bf abstract weak element} $[L]$ and {\bf abstract element} $\langle L \rangle$ for an element $L \in \mathcal{F}$:   
\[
  [L] := \begin{cases}
    \text{the inverse image by } q \\
    \text{of the connected component of} \\
     \{ L' \in X/\mathcal{F} \mid L \cong L', \overline{L} - L = \overline{L'} - L' \} 
     \text{ containing } L & \text{if } L \subseteq \mathrm{Cl}(\mathcal{F}) \sqcup \mathrm{P}(\mathcal{F}) \\
    \text{the inverse image by } q \\
    \text{of the connected component of} \\
  \{ L' \in X/\mathcal{F} \mid L \cong L', \overline{L} = \overline{L'} \}  \text{ containing } L & \text{if } L \subseteq \mathrm{R}(\mathcal{F}) 
  \end{cases}
\]
\[
  \langle L \rangle := \begin{cases}
    \text{the inverse image by } q \\
    \text{of the connected component of} \\
     \{ L' \in X/\mathcal{F} \mid  \overline{L} - L = \overline{L'} - L' \} 
     \text{ containing } L & \text{if } L \subseteq \mathrm{Cl}(\mathcal{F}) \sqcup \mathrm{P}(\mathcal{F}) \\
    \text{the inverse image by } q \\
    \text{of the connected component of} \\
    \{ L' \in X/\mathcal{F} \mid \overline{L} = \overline{L'} \} \text{ containing } L & \text{if } L \subseteq \mathrm{R}(\mathcal{F}) 
  \end{cases}
\]
Here $ L \cong L'$ means that $L$ and $L'$ are homeomorphic. 
\end{definition}

For a point $x \in X$, define $[x] := [\F(x)]$ and $\langle x \rangle := \langle \F(x) \rangle$.
As the same argument of the proof of Lemma~\ref{lem:ch_ab_ele}, we have the following observation.

\begin{lemma}\label{lem:ch_ab_ele_decop02}
The abstract elements form an invariant decomposition and satisfy the following property: 
\[
  \langle x \rangle = \begin{cases}
        \text{the inverse image by } q \text{ of the connected component of} \\
   \mathop{\mathrm{Cl}}(\F)/\F \text{ containing } \F(x) & \text{if } x \in\mathop{\mathrm{Cl}}(\F) \\
      \text{the inverse image by } q \text{ of the connected component of} \\
   q(\{ x' \in \mathrm{P}(\mathcal{F}) \mid \overline{\F(x)} - \F(x) = \overline{\F(x')} - \F(x') \}) \\
     \text{containing } \F(x) & \text{if }  \in \mathrm{P}(\mathcal{F})\\
     \text{the inverse image by } q \text{ of the connected component of} \\
  q(\{ x' \in \mathrm{R}(\mathcal{F}) \mid \overline{\F(x)} = \overline{\F(x')} \})  \text{ containing } \F(x) & \text{if } x \in \mathrm{R}(\mathcal{F}) 
  \end{cases}
\]
\end{lemma}

\begin{proof}
Let $\mathcal{F}$ be an invariant decomposition on a topological space $X$. 
Fix a point $x \in X$. 
If $x \in\mathop{\mathrm{Cl}}(\F)$, then $\langle x \rangle$ is the inverse image by $q$ of the connected component of $q(\{ x' \in X \mid \overline{\F(x)} - \F(x) = \overline{\F(x')} - \F(x') = \emptyset \}) = q(\mathop{\mathrm{Cl}}(\F)) = \mathop{\mathrm{Cl}}(\F)/\F$. 

Suppose that $x \in \mathrm{R}(\F)$. 
If there is a point $y \in \langle x \rangle \cap\mathop{\mathrm{Cl}}(\F)$, then $x \in \overline{\F(y)} = \F(y) \subseteq\mathop{\mathrm{Cl}}(\F)$, which contradicts $x \in \mathrm{R}(\F) = X - (\mathrm{Cl}(\F) \sqcup \mathrm{P}(\F))$. 
Thus $\langle x \rangle \cap\mathop{\mathrm{Cl}}(\F) = \emptyset$. 
We claim that $\langle x \rangle \cap \mathrm{P}(\F) = \emptyset$. 
Indeed, assume that there is a point $y \in \langle x \rangle \cap \mathrm{P}(\F)$. 
Since $x \in \mathrm{R}(\F)$, we have $\F(x) \subseteq \overline{\F(y)} - \F(y)$. 
Since $\overline{\F(y)} - \F(y)$ is closed, we have $\overline{\F(y)} = \overline{\F(x)} \subseteq \overline{\F(y)} - \F(y)$, which is a contradiction. 
Thus $\langle x \rangle$ is the inverse image by $q$ of the connected component of $q(\{ x' \in \mathrm{R}(\F) \mid \overline{\F(x)}  = \overline{\F(x')} \})$. 

Suppose that $x \in \mathrm{P}(\F)$. 
Then the derived set $\overline{\F(x)} - \F(x) \neq \emptyset$ is closed, and the abstract element $\langle x \rangle$ is the inverse image by $q$ of the connected component of $q(\{ x' \in X \mid \overline{\F(x)} - \F(x) = \overline{\F(x')} - \F(x') \neq \emptyset \})$. 
Therefore $\langle x \rangle \cap\mathop{\mathrm{Cl}}(\F) = \emptyset$. 
We claim that $\langle x \rangle \cap \mathrm{R}(\F) = \emptyset$. 
Indeed, assume that there is a point $y \in \langle x \rangle \cap \mathrm{R}(\F)$. 
Since $\overline{\F(x)} - \F(x)$ is closed, the set difference $\overline{\F(y)} - \F(y) = \overline{\F(x)} - \F(x)$ is closed and so $y \in \mathrm{P}(\F)$, which contradicts $y \in \mathrm{R}(\F)$. 
Thus $\langle x \rangle$ is the inverse image by $q$ of the connected component of $q(\{ x' \in \mathrm{P}(\F) \mid \overline{\F(x)} - \F(x) = \overline{\F(x')} - \F(x') \})$. 
\end{proof}

The previous lemma implies the following statement.

\begin{lemma}\label{lem:ch_ab_ele_decop}
The abstract weak elements form an invariant decomposition and satisfy the following property: 
\[
  [x] = \begin{cases}
     \text{the inverse image by } q \text{ of the connected component of} \\
   q(\{ x' \in \mathrm{Cl}(\mathcal{F}) \mid \F(x) \cong \F(x') \}) 
      \text{ containing } \F(x) & \text{if } x \in \mathrm{Cl}(\mathcal{F}) \\
     \text{the inverse image by } q \text{ of the connected component of} \\
    q(\{ x' \in \mathrm{P}(\mathcal{F}) \mid \F(x) \cong \F(x'), \overline{\F(x)} - \F(x) = \overline{\F(x')} - \F(x') \}) \hspace{-7pt} {} \\
     \text{containing } \F(x) & \text{if } x \in \mathrm{P}(\mathcal{F})\\
     \text{the inverse image by } q \text{ of the connected component of} \\
  q(\{ x' \in \mathrm{R}(\mathcal{F}) \mid  \F(x) \cong \F(x'), \overline{\F(x)} = \overline{\F(x')} \})  \\\text{containing } \F(x) & \text{if } x \in \mathrm{R}(\mathcal{F}) 
  \end{cases}
\]
\end{lemma}

Define the {\bf abstract weak element space} $X/[\F]$ as a quotient space $X/\sim_{[\F]}$ defined by $x \sim_{[\F]} y$ if $[x] = [y]$.
Similarly, define the {\bf abstract element space} $X/\langle \F \rangle$ as a quotient space $X/\sim_{\langle \F \rangle}$ defined by $x \sim_{\langle \F \rangle} y$ if $\langle x \rangle = \langle y \rangle$.
Since $\F(x) \subseteq [x] \subseteq \langle x \rangle$ for any $x \in X$, the abstract weak element space is a quotient space of the decomposition space, and the abstract element space is a quotient space of the abstract weak element space.

\subsubsection{Quasi-recurrence of elements of a decomposition}

A point $x \in X$ is {\bf non-maximal} if there is a point $y \in X$ such that $\overline{\F(x)} \subsetneq \overline{\F(y)}$. 
Define quasi-recurrence as follows. 

\begin{definition}
A point $x$ of a topological space $X$ with a decomposition $\F$ is {\bf quasi-recurrent} if $\langle x \rangle$ contains either a point in $\mathcal{R} (\F)$ or a non-maximal point. 
\end{definition}

Denote by $\bm{\mathop{Q} (\F)}$ the set of quasi-recurrent points, called the {\bf quasi-recurrent set} of a decomposition $\F$ and by $\max \F$ the set of maximal point with respect to the saturated topology. 
As the same argument of the proof of Lemma~\ref{lem:ch_fol_rec_chain_top}, we have the following statement.

\begin{lemma}\label{lem:ch_fol_rec_chain_decomp}
The following statements hold for an invariant decomposition $\F$ on a topological space $X$: 
\\
{\rm(1)} $\mathcal{R} (\F) \sqcup (\bigcup_{L \in \F} (\overline{L} - L) \cap \mathrm{P}(\F)) \subseteq \mathcal{Q}(\F)$. 
\\
{\rm(2)} $X - \mathcal{Q}(\F) \subseteq \mathrm{P}(\F) \cap \max \F$. 
\\
{\rm(3)} The set $\mathcal{Q}(\F)$ of quasi-recurrent points and the set of non-maximal points are $\hat{\F}$-invariant. 
\\
{\rm(4)} For any element $L \subset \mathrm{P}(\F) \cap \max \F$, the derived subset $\overline{L} - L$ is a closed $\hat{\F}$-invariant subset contained in $\mathcal{Q}(\F)$. 
\end{lemma}

\begin{proof}
Let $\F$ be an invariant decomposition on a topological space $X$. 
For any non-maximal point $x \in X$, there is a point $y \in X$ such that $\hat{\F}(x) \subseteq \overline{\F(x)} \subsetneq \overline{\F(y)}$. 
This means that the set of non-maximal points is $\hat{\F}$-invariant. 
By definition of quasi-recurrence, from Lemma~\ref{lem:saturated}, the set $\mathcal{Q}(\F)$ of quasi-recurrent points is $\hat{\F}$-invariant such that $\bigcup_{L \in \F} (\overline{L} - L) \subseteq \mathcal{Q}(\F)$, $X - \mathcal{Q}(\F) \subseteq \max \F$, and $X - \mathrm{P}(\F) = \mathcal{R} (\F) \subseteq \mathcal{Q}(\F)$. 
Then $\mathcal{R} (\F) \sqcup (\bigcup_{L \in \F} (\overline{L} - L) \cap \mathrm{P}(\F)) \subseteq \mathcal{Q}(\F)$ and $X - \mathcal{Q}(\F) \subseteq \mathrm{P}(\F) \cap \max \F$.
 
For any element $L \subset \mathrm{P}(\F) \cap \max \F$, since any closed $\F$-invariant subset is $\hat{\F}$-invariant,  the derived subset $\overline{L} - L$ is a closed $\hat{\F}$-invariant subset in $X - \max \F \subseteq \mathcal{Q}(\F)$. 
\end{proof}

The equality in Lemma~\ref{lem:ch_fol_rec_chain_decomp} {\rm(1)} does not hold in general (see Example in \S \ref{ex:04}).

\subsubsection{Morse hyper-graph of a  decomposition}

As the case of topological spaces, we define the Morse hyper-graph of a decomposition on a topological space as follows: 
Let $\F$ be a decomposition on a topological space $X$ with a set $\mathcal{M} = \{ M_\lambda \}_{\lambda \in \Lambda}$ of disjoint $\hat{\F}$-invariant non-empty subsets $M_{\lambda}$ ($\lambda \in \Lambda$). 
For any $I \subseteq \lambda$, define a {\bf hyper-edge} $H_{I}$ as follows: For any point $x \in X - \bigsqcup_{\lambda \in \Lambda} M_\lambda$, we say that $x \in H_{I}$ if there are disjoint $\hat{\F}$-invariant non-empty subsets $(C_i)_{i \in I}$ of the subset $\overline{\F(x)} - \F(x)$ such that $C_i \subseteq M_i$ and $\overline{\F(x)} - \F(x) = \bigsqcup_{i \in I} C_i$. 
Put $V := \{ M_{\lambda} \mid \lambda \in \Lambda \}$ and $H := \{ \{ M_i \}_{i \in I} \mid H_{I} \neq \emptyset, I \subseteq \Lambda \}$. 

\begin{definition}
A hyper-graph $\mathcal{G}_{\mathcal{M}} := (V, H)$ is the {\bf Morse hyper-graph} of $\mathcal{M}$ if $M - \bigsqcup_{\lambda \in \Lambda} M_\lambda = \bigsqcup_{I \subseteq \Lambda} H_{I}$. 
\end{definition}

For a Morse hyper-graph $\mathcal{G}_{\mathcal{M}} = (V, H)$, the associated graph $G_{\mathcal{M}} = (V, E)$ is defined by $E:= \{ \{ M_i, M_j \} \mid H_{I} \neq \emptyset, i, j \in I \subseteq \Lambda \}$. 
We define the Morse hyper-graph of a decomposition as follows. 

\begin{definition}
The hyper-graph $\mathcal{G}_{\mathcal{M}}$ is the {\bf Morse hyper-graph} of $\F$ if $\mathcal{M}$ is the set of inverse images by $q$ of connected components of $q(\mathcal{Q}(\F)) \subset X/\F$, where $q \colon X \to X/\F$ is the quotient map. 
\end{definition}
Then denoted by $\mathcal{G}_{\F}$ the Morse hyper-graph of $\F$. 
As the same argument of the proof of Lemma~\ref{lem:ex_MH_top_02}, the following existence of Morse hyper-graph of the decomposition holds.

\begin{lemma}\label{lem:ex_MH_top_03}
Let be $\F$ an invariant decomposition on a topological space $X$, $\mathcal{Q}' \subseteq X$ a $\hat{\F}$-invariant subset containing $\mathcal{Q}(\F)$, and $\mathcal{M} = \{ M_\lambda \}_{\lambda \in \Lambda}$ the set of inverse images by $q$ of connected components of $\mathcal{Q}'/\F$. 
Then the Morse hyper-graph $\mathcal{G}_{\mathcal{M}}$ exists. 
\end{lemma}

\begin{proof}
Since any connected component of a subset is closed in the subset, any connected component $M_\lambda$ is $\hat{\F}$-invariant and so the family $\mathcal{M}$ is a set of disjoint $\hat{\F}$-invariant non-empty subsets. 
Lemma~\ref{lem:ch_fol_rec_chain_decomp} implies that $X - \bigsqcup_{\lambda \in \Lambda} M_\lambda = X - \mathcal{Q}' \subseteq X - \mathcal{Q}(\F) \subseteq  \mathrm{P}(\F) \cap \max \F$. 
Fix a point $x \in X - \bigsqcup_{\lambda \in \Lambda} M_\lambda \subseteq \mathrm{P}(\F) \cap \max \F$. 
By Lemma~\ref{lem:ch_fol_rec_chain_decomp}, the derived subset $\overline{\F(x)} - \F(x) \subseteq \mathcal{Q}(\tau) \subseteq \mathcal{Q}' =  \bigsqcup_{\lambda \in \Lambda} M_\lambda$ is closed $\hat{\F}$-invariant and so the intersections $C_\lambda := (\overline{\F(x)} - \F(x)) \cap M_\lambda$ are $\hat{\F}$-invariant.
Set $I := \{ \lambda \in \Lambda \mid C_\lambda \neq \emptyset \}$. 
Then $x \in  H_I$. 
This means that $X - \bigcup_{\lambda \in \Lambda} M_\lambda = X - \mathcal{Q}(\F) = \bigsqcup_{I \subseteq \Lambda} H_{I}$. 
\end{proof}

The previous lemma implies the following existence. 

\begin{corollary}\label{lem:ex_MH}
The Morse hyper-graph for an invariant decomposition on a topological space exists.
\end{corollary}

%

In general, the Morse hyper-graph $\mathcal{G}_{\mathcal{M}}$ even for $\mathcal{M}$ as in Lemma~\ref{lem:ex_MH_top_03} is not a quotient space of the abstract element space (see the example in \ref{ex:mhg}).

%

\subsection{Reductions to Morse hyper-graphs for invariant decompositions}

Theorem~\ref{th:Morse_reduction} implies the following reduction to Morse hyper-graphs for invariant decompositions on topological spaces. 

\begin{theorem}\label{th:Morse_reduction_decomp}
The Morse hyper-graph for an invariant decomposition on a topological space is quotient spaces of the abstract element space and the abstract weak element space.
\end{theorem}

\begin{proof}
Let $\mathcal{G}_{\F} = (V, H)$ be the Morse hyper-graph for an invariant decomposition $\F$ with $V = \{ M_\lambda \}_{\lambda \in \Lambda}$ and $H = \{ \{ M_i \}_{i \in I} \mid H_{I} \neq \emptyset, I \subseteq \Lambda \}$.
Then $\mathcal{M} = \{ M_\lambda \}_{\lambda \in \Lambda}$ is the set of inverse images by $q$ of connected components of $q(\mathcal{Q}(\F)) \subset X/\F$.  
By $\mathcal{R} (\F) \subseteq \mathcal{Q}(\F) = \bigsqcup_{\lambda \in \Lambda} M_\lambda$, since the images by the quotient map $q \colon X \to X/\F$ of any abstract elements are connected, any recurrent abstract elements are contained in some $M_\lambda$. 

Fix any point $x \in \mathrm{P}(\F)$. 
Then the abstract element $\langle x \rangle$ is the inverse image by $q$ of a connected component of $q(\{ x' \in \mathrm{P}(\F) \mid \overline{\F(x)} - \F(x) = \overline{\F(x')} - \F(x') \neq \emptyset \})$ and the derived set $\overline{\F(x)} - \F(x)$ is closed and $\hat{\F}$-invariant. 
Suppose that $x$ is quasi-recurrent. 
Then $\langle x \rangle$ is contained in $\mathcal{Q}(\F) = \bigsqcup_{\lambda \in \Lambda} M_\lambda$ such that $q(\langle x \rangle)$ is connected in $\mathcal{Q}(\F)/\F$. 
Therefore there is an index $i \in \Lambda$ such that $\langle x \rangle \subseteq M_i$. 
Suppose that $x$ is not quasi-recurrent. 
Then $\langle x \rangle \subseteq X - \mathcal{Q}(\F) = X - \bigsqcup_{\lambda \in \Lambda} M_\lambda  = \bigsqcup_{I \subseteq \Lambda} H_{I}$. 
Fix any $x' \in \langle x \rangle$. 
Then $\overline{\F(x)} - \F(x) = \overline{\F(x')} - \F(x') \subseteq \mathcal{Q}(\F) = \bigsqcup_{\lambda \in \Lambda} M_\lambda$. 
Put $C_i := M_i \cap (\overline{\F(x)} - \F(x)) = M_i \cap (\overline{\F(x')} - \F(x'))$. 
Define $I := \{ i \in \Lambda \mid C_i \neq \emptyset \}$. 
Since $M_i$ and $\overline{\F(x)} - \F(x)$ are $\hat{\F}$-invariant, so is the intersection $C_i$ for any $i \in I$. 
Moreover, we have that $C_i \subseteq M_i$ and $\overline{\F(x')} - \F(x') = \overline{\F(x)} - \F(x) =  (\overline{\F(x)} - \F(x)) \cap \bigsqcup_{\lambda \in \Lambda} M_\lambda = \bigsqcup_{i \in I} C_i$. 
Then $x', x \in H_I$ and so $ \langle x \rangle \subseteq H_I$. 
\end{proof}

\subsection{Generalization of Reeb graphs of Morse functions}

We recall Reeb graph as follows. 
For a function $f \colon  X \to \R$ on a topological space $X$, the {\bf Reeb graph} of a function $f \colon  X \to \R$ on $X$ is a quotient space $X/\sim_{\mathrm{Reeb}}$ defined by $x \sim_{\mathrm{Reeb}} y$ if there are a number $c \in \R$ and a connected component of $f^{-1}(c)$ which contains $x$ and $y$.
Then the inverse image of a value of $\R$ is called the {\bf level set}. 
Notice that the Reeb graph of a Morse function (or more generally a function with finitely many critical points) on a closed manifold is a finite graph (see \cite[Theorem~3.1]{saeki2020reeb} for details). 
Moreover, by the proof of Proposition~\ref{prop:reeb}, note that any edge of the Reeb graph of a Morse function on a closed manifold is the leaf space of a codimension one product compact foliation, because the leaf space of a continuous codimension two (and so one) compact foliation of a compact manifold is Hausdorff \cite{EMS,epstein1972per,Epstein1976,Vogt1976}. 
We show that the abstract element space is a natural generalization of the Reeb graph of a Morse function. 

\begin{proposition}\label{prop:reeb}
The Reeb graph of a Morse function on a closed manifold is the abstract weak element space of the set of connected components of level sets as abstract multi-graphs. 
\end{proposition}

\begin{proof}
Let $f$ be a $C^1$ function with finitely many critical points on a closed manifold $M$ and $\F$ the set of connected components of level sets of $f$. 
Then any elements of $\F$ are closed. 
Denote by $\mathrm{Crit}(f)$ the set of critical points.
We claim that any elements of $\F$ contained in a connected component of the complement $X - \bigcup_{x \in \mathrm{Crit}(f)} \F(x)$ are homeomorphic to each other and are submanifolds. 
Indeed, fix such a connected component $C$.  
Then the restriction $f\arrowvert_C$ is a submersion, and so the restriction $\F\arrowvert_C$ is a codimension one foliation on the connected manifold $C$. 
Any element contained in $C$ of $\F$ is a closed codimension one submanifold because of the inverse function theorem. 
By the existence of the function $f$, the codimension one compact foliation $\F\arrowvert_C$ is transversely orientable, and so the holonomy group of each leaf of $\F\arrowvert_C$ is trivial.  
The Reeb stability (cf. \cite[Theorem~2]{thurston1974generalization}) implies that any leaf $L \in \F\arrowvert_C$ has its \nbd on which the restriction of $\F$ is a product foliation. 
Thus any leaves of $\F\arrowvert_C$ are homeomorphic to each other. 

Therefore any connected component of the complement $X - \bigcup_{x \in \mathrm{Crit}(f)} \F(x)$ is contained in the abstract weak element. 
Suppose that $f$ is Morse.
By the Morse lemma, any element of $\F$ containing a critical point is isolated but is not a manifold, and so such elements are abstract weak elements. 
This implies that any connected components of the complement $X - \bigcup_{x \in \mathrm{Crit}(f)} \F(x)$ correspond to abstract weak elements. 
Therefore any vertices and edges of the Reeb graph are abstract weak elements. 
This means that the Reeb graph of $f$ is the abstract weak element space $M/\F$ as an abstract multi-graph. 
\end{proof}



\subsection{Foliated spaces}
Recall some concepts of foliated spaces. 
Fix any $r \in\mathbb{Z}_{\geq 0} \sqcup \{ \infty \}$ and any $n \in\mathbb{Z}_{\geq 0}$. 
Let $Z, Z'$ be topological spaces and $U$ an open subset in $\R^n\times Z$ with coordinates $(x,z)$.
A mapping $f \colon U\to \R^n$ is of $\bm{C^r}$ if its partial derivatives up to order $r$ with respect to $x$ exist and are continuous on $U$. 
A mapping $h = (h_1, h_2) \colon U\to\R^n\times Z'$ of the form $h(x,z)=(h_1(x,z),h_2(z))$ is of $C^r$ if $h_1$ is of class $\bm{C^r}$ and $h_2$ is continuous, where .
\begin{definition}
A collection $\mathcal{U}=\{(U_i,\phi_i)\}_{i \in \Lambda}$ if a {\bf $C^r$ foliated atlas} of dimension $n$ on a topological space $X$ if it satisfies the following conditions: 
\\
{\rm(1)} The subset $\{U_i\}_{i \in \Lambda}$ is an open covering of $X$. 
\\
{\rm(2)} For any $i \in \Lambda$, there are a locally compact separable completely metrizable space $Z_i$ and an open ball $B_i$ in $\R^n$ such that $\phi_i \colon U_i\to B_i\times Z_i$ is a homeomorphism. 
\\
{\rm(3)} For any $i,j \in \Lambda$,  the coordinate changes $\phi_j \circ \phi_i^{-1}\arrowvert :\phi_i(U_i\cap U_j)\to\phi_j(U_i\cap U_j)$ are locally $C^r$ mappings of the form
\[
\phi_j \circ \phi_i^{-1}(x,z) = (g_{ij}(x,z),h_{ij}(z))\;.
\]
\end{definition}

Then $h_{ij}$ is called the {\bf local transverse components} of the changes of coordinates. 
Each $(U_i,\phi_i)$ is called a {\bf $C^r$ foliated chart}, and the inverse images $\phi_i^{-1}(B_i\times \{z\})$ for any $z\in Z_i$ are called {\bf plaques}. 
The foliated atlas of a topological space $X$ induces a locally Euclidean topology $\tau_{\mathcal{U}}$ on $X$. 
In fact, the basic open subsets are the plaques of all foliated charts. 
The connected components of $X$ with respect to this topology $\tau_{\mathcal{U}}$ are called {\bf leaves}. 
Each leaf is a connected $n$-manifold with the $C^r$ structure canonically induced by $\mathcal{U}$. 
The set of leaves is denoted by $\F$ and called a {\bf $C^r$ foliated structure} of dimension $n$ on $X$. 
The quotient space is called the {\bf leaf space} of $\F$ and denoted by $X/\mathcal{F}$. 
Notice that the leaf space is a decomposition space. 
For a point $x \in X$, denote by $\F(x)$ the leaf containing $x$. 
\begin{definition}
Two $C^r$ foliated atlases on $X$ define the same $C^r$ foliated structure if their union is a $C^r$ foliated atlas. 
\end{definition}

\begin{definition}
A maximal $C^r$ foliated atlas is called a {\bf $C^r$ foliated structure}.
\end{definition}

\begin{definition}
The pair of a topological space $X$ and its $C^r$ foliated structure is called a {\bf foliated space}.
\end{definition}
 
The $C^r$ foliated space with boundary is defined similarly, and the boundary of a $C^r$ foliated space is a $C^r$ foliated space without boundary.
A subset of a foliated space is {\bf saturated} (or {\bf invariant}) if it is a union of leaves. 
Notice that any saturated subspace is also a $C^r$ foliated space.

\subsubsection{Existence of Morse hyper-graphs of foliated spaces}

We have the following observations. 

\begin{lemma}\label{prop:02}
If $\F$ is either a foliated structure of a foliated space or the set of orbits of a group-action on a topological space, then $\F$ is an invariant decomposition. 
\end{lemma}

\begin{proof}
Since the closure of the saturation of any subset is $\F$-invariant, Lemma~\ref{lem:001} and Lemma~\ref{cor:correspondence} imply the assertion. 
\end{proof}

We have the following existence of Morse hyper-graphs of foliated spaces. 

\begin{theorem}\label{prop:Morse_fol}
Let $(X, \F)$ be either a compact foliated space or a topological space equipped with the set of orbits of a group action.  
The following statements hold: 
\\
{\rm(1)} The Morse hyper-graph $\mathcal{G}_{\F}$ exists.
\\
{\rm(2)} The Morse hyper-graph $\mathcal{G}_{\F}$ is the quotient spaces of the abstract element space $X/\langle \F \rangle$ and the abstract weak element space $X/[\F]$.
\end{theorem}

\begin{proof}
Lemma~\ref{prop:02} implies that $\F$ is invariant. 
By Corollary~\ref{lem:ex_MH}, the Morse hyper-graph $\mathcal{G}_{\F}$ exists. 
Theorem~\ref{th:Morse_reduction_decomp} implies assertions {\rm(2)}. 
\end{proof}

%

\section{Examples}

We state the following recurrent properties. 

\subsection{Recurrent properties}\label{sec:01a}

\subsubsection{Morse hyper-graphs of $\mathcal{M}$ which are not quotient spaces of abstract element space}\label{ex:mhg}

\begin{figure}
\begin{center}
\includegraphics[scale=0.125]{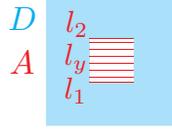}
\end{center}
\caption{A decomposition on a square $D = [0, 3]^2$}
\label{fig:quotient_mhg}
\end{figure}

There is the Morse hyper-graph $\mathcal{G}_{\mathcal{M}}$ for $\mathcal{M}$ is not a quotient space of the abstract element space. 
Indeed, let $D := [0,3]^2$ be a closed square, $A := (1,2) \times [1,2]$ a square, and $\F := \{ \{x \} \mid x \in D - A \} \sqcup \{ l_y \mid y \in [1,2] \}$ a decomposition as in Figure~\ref{fig:quotient_mhg}, where $l_y := (1,2) \times \{y\}$ are open intervals.
Then $\mathcal{Q}(\F) = \mathop{\mathrm{Cl}}(\F) = D -A$, $\mathrm{P}(\F) = A$, and $\overline{\mathcal{Q}(\F)} - \mathcal{Q}(\F) = A - \mathrm{int} A = l_1 \sqcup l_2$. 
Therefore the abstract (weak) element space is $D/[\F] = D/\langle \F \rangle = \{ D - A, A\}$ and the Morse hyper-graph $\mathcal{G}_{\F}$ of $\F$ is $(V, H) = (\{ M_1 = D-A\}, \{H_1 = A\})$. 
On the other hand, the family $\mathcal{M} := \{ D - (1,2)^2 \}$ is  the set of connected components of $\overline{\mathcal{Q}(\F)} = D - \mathrm{int} A = D - (1,2)^2$. 
Then the Morse hyper-graph $\mathcal{G}_{\mathcal{M}}$ for $\mathcal{M}$ is $(V_{\mathcal{M}}, H_{\mathcal{M}}) = (\{ M'_1 = D-\mathrm{int} A\}, \{H'_1 = \mathrm{int} A\})$. 
This means that the Morse hyper-graph $\mathcal{G}_{\mathcal{M}}$ is not a quotient space of the abstract element space $D/[\F] = \{ D - A, A\}$. 

Denote by $(X, \tau)$ the element space $D/\F = \{ \{ x \} \mid x \in D - A\} \sqcup \{ l_y \mid y \in [1,2] \}$. 
Then $(X, \tau)$ is a topological space such that $\mathcal{Q}(\tau) = \mathop{\mathrm{Cl}}(\tau) = (D - A)/\F = \{ \{ x \} \mid x \in D - A\}$, $\mathrm{P}(\tau) = A/\F =  \{ l_y \mid y \in [1,2] \}$ and $\overline{\mathcal{Q}(\tau)} = (D - \mathrm{int} A)/\F = \{ \{ x \} \mid x \in D - A\} \sqcup \{ l_1, l_2 \}$. 
Moreover, we have $D/[\tau] = D/\langle \tau \rangle = \{ (D - A)/\F, A/\F\}$ and the Morse hyper-graph $\mathcal{G}_{\tau}$ of $\tau$ is $(V'', H'') = (\{ M''_1 = (D - A)/\F\}, \{H''_1 = A/\F\})$ and is isomorphic to $\mathcal{G}_{\F}$. 
On the other hand, the family $\mathcal{M}_\tau := \{ \overline{\mathcal{Q}(\tau)} \}$ the set of connected components of $\overline{\mathcal{Q}(\tau)}$. 
Then the Morse hyper-graph $\mathcal{G}_{\mathcal{M}_\tau}$ for $\mathcal{M}_\tau$ is $(V_{\mathcal{M}_\tau}, H_{\mathcal{M}_\tau}) = (\{ M'''_1 = (D-\mathrm{int} A)/\F\}, \{H'''_1 = \mathrm{int} A/\F \})$ and is isomorphic to $\mathcal{G}_{\mathcal{M}}$. 
This means that the Morse hyper-graph $\mathcal{G}_{\mathcal{M}_\tau}$ is not a quotient space of the abstract element space $X/[\tau] = \{ (D - A)/\F, A/\F\}$.

\subsubsection{Non-closedness of quasi-recurrent set}\label{sec:01}

There is a topological space whose quasi-recurrent set is not closed. 
In fact, define a poset $(P = \Z \times \{ 0,1 \}, \leq)$ by $(n_1, n_2) < (m_1, m_2)$ if $n_1 = m_1$, $n_2 = 0$, and $m_2 = 1$. 
The family $\tau := \{ \emptyset \} \sqcup \{ P - F \mid F \text{ is a finite downset} \}$ is a topology on $P$ whose specialization pre-order is the partial order $\leq$. 
Then $\mathrm{P}(\tau) = \max \tau = \Z \times \{1\}$ is infinite and so is not open. 
Therefore $\Z \times \{0\} = \mathrm{Cl}(\tau) = \mathcal{Q}(\tau)$ is not closed,

\subsection{Necessary conditions}\label{sec:02}

We construct some examples to describe the necessity of conditions. 

\subsubsection{Necessity of invariance to become topologies}\label{ex:01}
There is a decomposition whose set of $\F$-invariant open subsets is not a topology. 
In fact, for a decomposition $\F := \{ \{0 \} \times [0,1] \} \sqcup \{ (x, y) \mid x \neq 0, y \in [0,1] \}$ on a closed disk $[0,1]^2$, the set $\tau_{\F}$ is not a topology.

\subsubsection{Necessity of local compactness for correspondence of recurrences}\label{ex:02}
There is a transitive flow $v$ on a metrizable space with non-$\tau_{X/v}$-recurrent but $v$-recurrent points. 
Indeed, applying a dump function to an irrational rotation, consider a vector field on $\T^2$ with one singular point $x$ such that each non-singular orbit is dense. 
Let $X  := O \sqcup \{ x \}$ be the union of $x$ and a non-singular orbit $O$ and $v$ the restriction of the flow. 
Then $X$ is metrizable and $v$ consists of one $v$-recurrent non-closed orbit $O$ and one singular point $x$. 
Thus $\overline{O} -O = \{ x \}$ is closed and so $O$ is not $\tau_{X/v}$-recurrent but $v$-recurrent.  

\subsubsection{Necessity for using abstract elements in the definition of quasi-recurrence}\label{ex:04} 

\begin{figure}
\begin{center}
\includegraphics[scale=0.2]{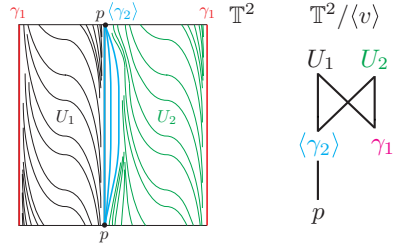}
\end{center}
\caption{A decomposition on a torus and its abstract element space.}
\label{fig:toral_flow}
\end{figure}

There is a singular codimension one foliation on a torus $\mathbb{T}^2$, which can be generated by a flow, with $\mathbb{T}^2 - \mathrm{P}(\F) \subsetneq \mathcal{R} (\F) \sqcup (\bigcup_{L \in \F} (\overline{L} - L) \cap \mathrm{P}(\F)) \subsetneq \mathbb{T}^2 -(U_1 \sqcup U_2) = \gamma_1 \sqcup \{ p \} \sqcup \langle \gamma_2 \rangle = \mathcal{Q}(\F)$ as in Figure~\ref{fig:toral_flow}. 
Put $X := \mathbb{T}^2/\F$. 
Then $X - \mathrm{P}(\tau_\F) \subsetneq \mathcal{R} (\tau_\F) \sqcup (\bigcup_{L \in X} (\overline{\{ L\}} - \{ L\}) \cap \mathrm{P}(\tau_\F)) \subsetneq X -(U_1 \sqcup U_2)/\F = \{ \gamma_1, \{ p \}, \langle \gamma_2 \rangle \} = \mathcal{Q}(\tau_\F)$. 
Note that the abstract element space $\mathbb{T}^2/\langle \F \rangle$ is weakly homotopic to a circle. 


\bibliographystyle{abbrv}
\bibliography{yt20211011}

\end{document}